\documentclass{article}
\usepackage{amssymb, amsmath}
\usepackage[latin1]{inputenc}
\usepackage[affil-it]{authblk}
\usepackage{hyperref}
\usepackage[onehalfspacing]{setspace}
\usepackage{ulem}
\makeatletter
\def\imod#1{\allowbreak\mkern10mu({\operator@font mod}\,\,#1)}
\makeatother
\begin{document}
\title{Linear groups in Galois fields. \\ A case study of tacit circulation \\ of explicit knowledge}
\author{Fr{\'e}d{\'e}ric Brechenmacher%
\thanks{Electronic address: \texttt{frederic.brechenmacher@euler.univ-artois.fr} \\ \textit{Ce travail a bénéficié d'une aide de l'Agence Nationale de la Recherche : projet CaaF{\'E} (ANR-10-JCJC 0101)}}}
\affil{Univ. Lille Nord de France, F-59000 Lille \\ U. Artois, Laboratoire de math{\'e}matiques de Lens (EA 2462) \\ rue Jean Souvraz S.P. 18, F- 62300 Lens France}
 \date{}
\maketitle
\begin{abstract}
This preprint is the extended version of a paper that will be published in the proceedings of the Oberwolfach conference "Explicit vs tacit knowledge in mathematics" (January 2012).  It presents a case study on some algebraic researches at the turn of the twentieth century that involved mainly French and American authors. By investigating the collective dimensions of these works, this paper sheds light on the tension between the tacit and the explicit in the ways some groups of texts hold together, thereby constituting some shared algebraic cultures.\\
Although prominent algebraists such as Dickson made extensive references to papers published in France, and despite the roles played by algebra and arithmetic in the development of the American mathematical community, our knowledge of the circulations of knowledge between France and the United States at the beginning of the 20th century is still very limited. 
It is my aim to tackle such issues through the case study of a specific collective approach to finite group theory at the turn of the 20th century. This specific approach can  be understood as a shared algebraic culture based on the long run circulation of some specific procedures of decompositions of the analytic forms of substitutions. In this context, the general linear group was introduced as the maximal group in which  an elementary abelian group (i.e., the multiplicative group of a Galois field) is a normal subgroup. 
\end{abstract}
This paper aims at stressing some aspects of my works on the history of algebra in the 19th and 20th century, which are related to the  tacit vs the explicit in the ways some groups of texts hold together. These aspects raise issues as to the historian's choices of a corpus of reference and of a scale of analysis. They also address the more general problem of articulating the individual and collective dimensions of mathematics.\\ 
Although the category "algebra" points to some collective organizations of knowledge, this category took on changing identities in different times and spaces. Until the 1930s, "algebra" was especially not usually referring to an object-oriented discipline, i.e., as identifying both a corpus of specialized knowledge revolving around some specific objects and the institutionalized practices of transmissions of a group of professional specialists (the "algebraists").\cite{1} In France, for instance, algebra was, on the one hand, traditionally considered in the teaching of mathematics as an "elementary," or "intermediary," discipline encompassed by "the higher point of view" of analysis. On the other hand, algebra was also pointing to some procedures that made a "common link"  between researches in the various branches of the mathematical sciences.\cite{6} What was explicitly identified as  "algebraic" therefore often pointed to some implicit circulations between various theories. This situation makes it customary to study carefully the ways texts were referring one to another, thereby constituting some shared algebraic cultures.\\
I shall introduce this paper by making explicit how such issues came up in my research work, before focusing on a case study on "linear groups in Galois fields" at the turn of the 20th century. Although the latter designation may seem to make explicit some collective interests for a theory revolving around a specific object, i.e. $Gl_n(F_{p^n})$ ($p$ a prime number), we shall see that this designation actually supported the implicit reference to a specific algebraic culture that had developed over the course of the 19th century. 
\section{Introduction : the "versus"} Throughout the whole of 1874, Jordan and Kronecker were quarreling over the organization of the theory of bilinear forms. This controversy made explicit some tacit knowledge and know-how, such as some conflicting disciplinary ideals on the roles of algebra vs arithmetic, some epistemic values of simplicity vs effectivity, and two opposed internal philosophies of generality. \cite{2} But more importantly for the topic of the present paper, this quarrel also highlights the tacit intertextual relations that  lied beneath the expression "theory of bilinear forms." \\
The controversy  started with a public quarrel of priority over two theorems before turning into a private correspondence. As is illustrated by the following letter from Jordan to Kronecker (February 1874), the epistolary communication was mostly devoted to making explicit some tacit relations between some texts of Weierstrass, Christoffel, and Kronecker:
\begin{quote}
\sout{I would not like} ... \sout{that I want to wage war and} that I would prefer \sout{a controversy / war} some public debates to friendly discussions. I was not the one to \sout{open fire} start the controvery. It is true that I have published (as was my evident right) without consulting you some researches that were completing yours ... But if instead of going abruptly public with this issue, you would have contacted me for exchanging the friendly explanations I would have been legitimately expecting, we would certainly have reached an agreement. Following your indications, I would have \sout{immediately noticed, what I recognized too late, that your method of 1868} reread more carefully your memoir of 1868 and would have noticed, what I did not recognized at first sight, that bilinear forms were implicitly included in your work despite not being cited. \footnote{ \sout{Je ne voudrais pas} ..., \sout{que je désire la guerre}, et que je préférerais \sout{une polémique / guerre} des débats publics à des explications amicales. Ce n'est pas moi qui ait \sout{ouvert les hostilités} commencé la polémique. J'ai publié il est vrai (c'était mon droit évident) sans vous consulter des recherches qui complétaient les vôtres... Si au lieu de jeter brusquement ce débat dans le public, vous vous étiez adressé à moi pour échanger des explications, comme je me voyais en droit de l'espérer, nous nous serions sans doute entendu. Sur votre indication, j'aurais \sout{constaté immédiatement, ce que j'ai  reconnu trop tard, que votre méthode de 1868} relu plus attentivement votre mémoire de 1868 et constaté, ce que je n?avais pas remarqué à première vue, que les formes bilinéaires non citées dans votre travail, y sont pourtant implicitement comprises.} 
\end{quote}
 After Jordan had made himself familiar with the collective of texts Kronecker made explicit, the quarrel  eventually went public once again. These successive episodes highlight that even though the "theory of bilinear forms" may have been considered at first sight as the explicit reference to  a collective of mathematical methods and notions, it was through some tacit intertextual relations that this theory was making sense.  Moreover, the tacit dimension of the theory was recognized by the actors themselves. Both Jordan and Kronecker indeed agreed that it was no more than a "slight fault" ("levis culpa") to fail to grasp the relevant underlying intertextual relations, which had thus to be made explicit through direct communication. Both protagonists also recognized having been already both responsible and victims of such slight faults in the past, as is illustrated by the following letter from Jordan to Kronecker (January 1874):
\begin{quote}
As I conclude this letter, please allow me Sir to express my regrets that you did not \sout{write to me your} ask me for some explanations before publishing your criticisms ... I especially regret your decision to involve in this matter the name of M. Christoffel ... I should add that I did not know his memoir at the time \sout{which is quite legitimate}. I am moreover convinced that [Christoffel] would have agreed with me that one has to be modest about such claims. It is by following this principle that I did not say anything against the publication in the last issue of M. Borchardt's journal of a voluminous paper by M. Sohnke on the symmetry in the plane, in which I was not cited even though I have completely tackled this issue for both the space and the plane five or six years ago in MM. Brioschi's and Cremona's journals. \footnote{ Permettez moi Monsieur, en terminant cette lettre, de vous exprimer le regret que vous ayez attendu pour \sout{m'écrire vos} demander des explications d'avoir publié vos critiques ... Je regrette surtout que vous ayez cru devoir introduire dans cette question le nom de M. Christoffel ... J'ajouterai qu'à cette époque, je ne connaissais pas son mémoire \sout{ce qui est bien permis}. Je suis persuadé d'ailleurs qu'il aurait été d'avis comme moi, qu'il faut être fort réservé dans les réclamations de ce genre. \sout{Aussi ai-je} C'est d'après ce principe que j'ai laissé passer sans rien dire, dans le dernier numéro du Journal de M. Borchardt, un mémoire volumineux de M. Sohnke sur la symétrique dans le plan où je n'étais pas cité, bien que j'eusse traité complètement cette question pour le plan et pour l'espace dans le journal de MM Brioschi et Cremona, il y a de cela cinq à six ans.} 
\end{quote}
\section{The tacit vs the explicit in building some networks of texts}
The 1874 controversy highlights the problem of the selection of the corpuses in which a given text is making sense, i.e., the identification  of some networks of texts. But such networks cannot be identified as webs of quotations. \cite{10} Not only do practices of quotations vary in times and spaces but intertextual relations may also be implicit. My approach to this problem consists in choosing a point of reference from which a first corpus is built by following  systematically the explicit traces of intertextual relations. A close reading of the texts of this corpus then gives access to some more implicit forms of references. For instance, in the case of the 1874 controversy, both protagonists were referring to the "equation to the secular inequalities in planetary theory" \cite{2}. Such an explicit reference implicitly pointed to  a network of texts that had been published over the course of the 19th century in various theoretical frameworks, e.g., celestial mechanics, analytical geometry, complex analysis, or arithmetic. The 1874 controversy opposed two attempts to turn this traditional algebraic culture into an object-oriented theory, which Jordan aimed to root on group theory while Kronecker laid the emphasis on the theory of quadratic forms.
\section{Linear groups in Galois fields}
In the framework of a collective research project, \footnote{ CaaFÉ : Circulations of algebraic and arithmetic practices and knowledge  (1870-1945) : France, Europe, U.S.A ; \href{ http://caafe.math.cnrs.fr}{http://caafe.math.cnrs.fr}} a database of intertextual references has been  worked out for all the texts published in algebra in France from 1870 to 1914. \footnote{ The corpus has been selected by using the classification of the \textit{Jahrbuch}. On Thamous database of intertextual references, see  \href{ http://thamous.univ-rennes1.fr/presentation.php}{http://thamous.univ-rennes1.fr/presentation.php}} One of the subgroups of this corpus gives rise to a coherent network of texts  which involved mainly French and American authors from 1893 to 1907. \footnote{ On the one hand, Jordan, Borel and Drach, Le Vavasseur, de Séguier, Autonne, etc. On the other hand, Moore, Dickson, Schottenfels, Wedderburn, Bussey, Börger,  Miller, Manning, etc.}  A quite representative example  is the book "{\'E}lements de la th{\'e}orie des groupes abstraits," published in 1904 by Jean Armand de S{\'e}guier, one of the main French authors of the group.\\
Let us  characterize further this network by looking at its main shared references. These were, on the one hand,  some French papers published in the 1860s, and, on the other hand, Moore's introduction of the abstract notion of Galois field in 1893 in addition to Dickson's 1901 monograph on \textit{Linear groups with an exposition of the Galois field theory}. We shall see that the two times and spaces involved here point to a shared algebraic culture that can neither be identified to a discipline nor to any simple national or institutional dimension.
\section{On nations and disciplines}
Moore's 1893 "abstract" notion of Galois field has been assumed to highlight the influence of "the abstract point of view then increasingly characteristic of trendsetting German mathematics" on the emergence of both the "Chicago algebraic research school" \cite{8} and the "American mathematical research community" \cite{9}. Here two kinds of categories have been used for making explicit some collective dimensions of mathematics, i.e., on the one hand, some national, or more local, collectives of mathematicians (the U.S.A., Germany, Chicago) and, on the other hand, some mathematical disciplines (abstract algebra). But even though the roles played by German universities in the training of many American mathematicians have been well documented, the influence of this institutional framework on mathematics has been assumed quite implicitly. Here two difficulties arise.\\ 
First, the role attributed to "abstract algebra" reflects the tacit assumption that the communication of some local tacit knowledge should require direct contact. Because the historiography of algebra has usually emphasized the abstract approaches developed in Germany, and especially in the center of G\"{o}ttingen, other, more local, abstract approaches, such as in Cambridge or Chicago, have raised issues about the imperfect communication of some tacit knowledge, as exemplified by the late interbreeding in the 1930s of the German \textit{Moderne Algebra} and the Anglo-American approach to associative algebras. \cite{11} In this framework, the algebraic developments in France, such as S{\'e}guier's, have been either ignored or considered as some isolated attempts modeled on German or Anglo-American approaches. \cite{12} \\
Second, both disciplines and nations are actors categories, which even though they were much involved in public discourses,  cannot usually be directly transposed to the collective dimensions of mathematical developments \cite{5}. Public discourses indeed often involved some boundary work that was not only setting delimitations between mathematicians and non mathematicians but was also supporting some hierarchies among the practitioners of mathematics (researchers vs teachers and engineers, pure vs applied, analysts vs algebraists, etc.). \footnote{ L.Turner has shown that "research" was a category Mittag Leffler recurrently appealed to in connection to the notions of "contribution" for establishing hierarchies not only between various stratas of practitioners of mathematics, or between a "general public" and a "special public" but also between nations themselves in Mittag Leffler's own construction of an international space \cite{15}}.  These boundaries often reflect the roles taken on by some authorities as well as the ways some individuals embodied some collective models of mathematical lives or persona. For instance, several authorities of mathematics in France recurrently expressed publicly an opposition between the unifying power of analysis and the poverty of considering algebra as an autonomous discipline, thereby blaming some approaches developed in Germany.  Picard was one of the main advocate of such an opposition. But even though he was publicly celebrating Jordan's approach to Galois theory, it was to Kronecker's approach to this theory that he was appealing in his own mathematical work. \\
Recall that Moore's 1893 paper was read at the congress that followed the World Columbian exposition in Chicago. The world fair was dedicated to the discovery of America and was the occasion of much display of national grandeur. \cite{13} In parallel to the elevation of the first great wheel, presented by the Americans as a challenge to the Eiffel tower, or to the viking ship that sailed from Norway as a counterpoint to the replica of Columbus's three caravels, the architectural influence of the French {\'E}cole des Beaux arts was challenged by the German folk village. The latter especially included an exhibit of the German universities in which Klein was delivering a series of lectures which aimed at "passing in review some of the principal phases of the most recent development of mathematical thought in Germany." Klein was also the glorious guest of the congress while Moore was  both the host of the congress and one of its main organizers. The latter's concluding lecture was a tribute to Klein's \textit{Icosahedron}. It indeed aimed at generalizing to a "new doubly infinite system of simple groups," i.e., $PSL_2(F_{p^n})$,  what was then designated as the three "Galois groups," i.e., $PSl_2(F_p),  p=5,7,11$, involved in the modular equations that had been investigated by Galois, Hermite, and Klein among others. \cite{9} The generalization consisted in having the analytic form of unimodular binary linear fractional substitutions $\frac{ax+b}{cx+d}$ operate on a finite "field" of letters indexed by Galois number theoretic imaginaries.\\
The nature of the relevant collective dimensions nevertheless change if one shifts the scale of analysis from institutions to texts. As we shall see, even though he had aimed at celebrating the emergence of some abstract researches in the U.S.A. in the framework of the G\"{o}ttingen tradition, Moore actually collided to the implicit collective dimension that was underlying the use of the analytic representation of substitutions on number theoretic imaginaries \cite{4}.\\
In 1893, Moore initially appealed to Klein-Fricke's 1890 short presentation of Galois imaginaries  but had not yet studied the references to Galois, Serret, Mathieu, Jordan, or Gierster which he cited from Klein-Fricke's texbook. On the one hand, what Moore designated as a Galois field actually corresponded to Cauchy's approach to higher congruences \cite{10}, as developed later by  Serret, and which Moore nevertheless attributed to Galois :
\begin{quote}
The most familiar instance of such a field ... is the system of $p$ incongruous classes (modulo $p$) of rational integral numbers. Galois discovered an important generalization of the preceding field ... the system of $p^n$ incongruous classes (modulo $p$, $F_n(x)$).
\end{quote}
 On the other hand, Moore's "abstract finite field" was actually close to Galois's approach. Moore's remark that "every finite field is in fact abstractly considered a Galois field" thus echoed the connection between two perspectives on number theoretic imaginaries, as it had  already been displayed in textbooks such as Serret's in 1866.  \cite{4} \\
 But even more dramatically, Moore's system of simple groups had actually already been introduced by Mathieu in 1861.  As a result, before the publication of the proceedings of the congress in 1896, Moore and his student Dickson struggled to access the tacit collective dimension of some texts published in France in the 1860s, especially by appealing to Jordan's 1870 \textit{Trait{\'e} des substitutions et des {\'e}quations alg{\'e}briques}. This appropriation resulted in the publication of a train of  papers on "Jordan's linear groups in Galois fields." \\
 Moore's 1893 paper thus eventually resulted in the circulation of some works that were foreign to Klein's legacy. This situation highlights the difficult problem of identifying the scales at which various forms of collective dimensions play a relevant role, especially in respect to the articulation of the collective dimensions of texts with the ones of actors, such as disciplines or nations. 
\section{The analytic representation of substitutions}
Let us now characterize more precisely the collective dimension to which Moore collided to in 1893. Given a substitution $S$ operating on $m$ letters $a_i$, the problem of the analytic representation consists in finding an analytic function $\phi$ such that $S(a_i)=\phi(i)$. Hermite's 1863 solution to this problem for the cases $m=5, 7$  would especially influence Dickson's 1896 thesis. But the analytic representation also requires an indexing of the letters. This issue had already been tackled in Poinsot's 1808 presentation of Gauss's decomposition of cyclotomic equations. \cite{10} As had been shown later by Galois in 1830, if $m=p^n$, the indices can be considered as the "imaginary solutions" of the congruence $x^{p^n}-x\equiv 0 \imod{p}$ that generalized the indexations given by the roots of Gauss's cyclotomic equations. \\
Such analytic forms were nevertheless not limited to representations  but were interlaced with some specific procedures. The procedure of reduction of  "linear" substitutions $(i, ai+b$) into combinations of cycles $(i, i+1)$ and $(i, gi)$ had especially played a key role in Galois's characterization of solvable irreducible equations of prime degree, which roots $(x_i, x_{\phi(i)})$ are permuted by substitutions of "a linear form" $(x_i, x_{ai+b})$ \cite{4}. In modern parlance, Galois's theorem and its proofs boil down to showing that the linear group is the maximal group in which an elementary abelian group (the cyclic group $F_p^*$ in the case $n=1$) is a normal subgroup.\footnote{ Given $(x_i, x_{\phi(i)})$  the analytic representation of the roots, then $\phi(i+a)=\phi(i)+A$ because the cyclic group of substitutions $(i, i+a)$ is a normal subgroup. Thus $\phi(i+2a)=\phi(i)+2A)$, ..., $\phi(i+ma)=\phi(i)+mA$. Thus $\phi(m)=Am+B$ is of a "linear form" (i.e., affine).}  It was in attempting to generalize this theorem to the analytic forms of roots $(x_{i, i', ..., i^{(n)}}, x_{\phi(i), \psi(i'), ..., \sigma(i^{(n)})})$ of equations of degree $p^n$ that Galois introduced the number theoretic imaginaries.\\
Later on, in the 1860s,  Jordan investigated further the group "originating" from the problem of finding the analytic form of the maximal group $T$ in which the group of substitutions $(i, i', ..., i^{(n)} ; i+a, i'+a', ... , i^{(n)}+a^{(n)})$ (i.e., $F_{p^n}^*$) is a normal subgroup. He followed Poinsot's reformulation of Gauss's decomposition of cyclotomic equations in dividing the letters on which the substitutions are acting into "groups" (i.e., blocks of imprimitivity, which  are sets on which the group is acting but that Jordan also called groups), each of a same cardinal $p^n$.  $T$ was then simultaneously partitioned into a "combination of displacements between the groups" [i.e. blocks] and of permutations of the letters within each of the groups [i.e., blocks]." As a result, the substitutions of $T$ were divided into two "species". On the one hand, inside each block, the letters were cyclically substituted by the first specie of substitutions $(i, gi)$ operating by the multiplication of a primitive root $g$ (mod.$p$). On the other hand, the second specie $(i, i+1)$ substituted cyclically the blocks themselves. Each specie of substitution corresponded to one of the two forms of representation of cycles. Their products generated linear forms $(i, ai+b)$ with $i \in F_{p^n}$, i.e., if $i=(i, i', i'',  ... , i^{(n)})$ :
\[
\begin{vmatrix}
\phi(i) \equiv  ai+bi'+ci''+... + vi^{(n)}+d \\
\phi(i') \equiv a'i+b'i'+c'i''+... + v'i^{(n)}+d' \\ 
\phi(i'') \equiv a''i+b''i'+c''i''+... + v''i^{(n)}+d'' \\ 
... ... \\ 
\end{vmatrix}
mod (p)
\]
$T$ was therefore a linear group, as  Jordan would designate it in the mid 1860. During the following decades, Jordan especially laid the emphasis on the procedures of reductions of the analytic representations of linear substitutions. The statement  of Jordan's canonical form theorem between 1868 and 1870 especially gave a generalization to the reduction of the analytic form $(i, ai+b)$ into $(i, i+1)$  and $(i, gi)$ by reducing any linear substitution on $F_{p^n}$ to the following form:
\[
\begin{vmatrix}
y_0, z_0, u_0, ..., y'_0, ... & K_0y_0, K_0(z_0+y_0), ... , K_0y'_0 \\
y_1, z_1, u_1, ..., y'_1, ... & K_1y_1, K_1(z_1+y_1), ... , K_1y'_1 \\
.... & ... \\
v_0, ... & K'_0v_0, ... \\
... & ... \\
\end{vmatrix}
\]
For Jordan, such reductions were the very "essence of the question." They were at the core of what Jordan, following Poinsot, was designating as the "theory of order,"\cite{16}  i.e., a theory dealing with the "relations" between classes of objects in contrast with concerns for specific objects. Reductions of analytic representations were indeed supporting links between various branches of mathematics such as number theory (cyclotomy),  the theory of substitutions, the "problem of the irrationals" (Galois theory), crystallography (symmetries), mechanics (rotation of a solid body), analysis situs (polyedrons, Riemann surfaces), differential equations, etc.\cite{6} The process of reduction of a linear group into chains of subgroups was for instance compared to the unscrewing of a helicoïdal motion into some rotation and translation motions on the model of Poinsot's interpretations of the two forms  $(i, gi)$ and $(i, i+1)$ of analytic representations of  cycles as two forms of motions on a circle.\\
Jordan's 1870 \textit{Trait{\'e}} played a key role in the development of a specific algebraic culture based on the reduction of the analytic representation of $n$-ary linear substitutions. This culture can not only be be traced in France in the works of authors such as Poincar{\'e}, Picard, Autonne, Cartan, S{\'e}guier, etc., \cite{Brechenmacher:2011b} but it also circulated in the U.S.A. after Dickson's 1896 thesis \cite{6}. 
\section{Conclusion}
The network of texts that revolved around "Jordan's linear groups in Galois field" at the turn of the 20th century had underlying it a specific algebraic culture based on procedures of reductions of the analytic forms of substitutions. It was because they shared this culture that some French and American authors were able to interact with each others, even though most of them did not have any direct contact and did not share any social framework, as is exemplified by such different figures as de S{\'e}guier, an aristocrat jesuit abbot, and Schottenfels, one of the first women to graduate in mathematics at Chicago. Communication was nevertheless partial and was actually mostly limited to some shared practices, such as the use of Jordan's canonical form. A telling example is the new formulation that was given repeatedly and indepently to Jordan's "origin" of the linear group as the theorem stating that the group of automorphisms of an elementary abelian group $ F_{p^n}^*$ is $Gl(F_{p^n })$ (Burnside, Moore, Levavasseur, Miller, Dickson, S{\'e}guier).\\
We have seen that the systematic investigation of explicit traces of intertextual relations also sheds light on some more implicit collective forms of references, such as the one that lied beneath expressions such as "linear groups in Galois fields." This situation highlights the crucial role played by some networks of texts in the shaping of some algebraic cultures at a time when "algebra" was not yet referring to an object-oriented discipline. 

\end{document}